\documentclass[10pt,twoside]{article}
\usepackage{graphicx}
\usepackage{amsmath}
\usepackage{amsfonts}  
\usepackage{Latex-document}

\title{\bf  Geometry of Symplectic Intersections\vskip 6mm}

\author{P. Biran\vspace*{-0.5cm}\thanks{School of Mathematical Sciences, Tel-Aviv University, Tel-Aviv 69978,
Israel. Email:\ \ \ \ \ \ \ \ \ \ \ \ \ \ \ \
biran@math.tau.ac.il}}
\date{\vspace{-8mm}}

\begin{document}
\maketitle

\thispagestyle{first} \setcounter{page}{241}

\begin{abstract}

\vskip 3mm

In this paper we survey several intersection and non-intersection phenomena appearing in the realm of symplectic
topology. We discuss their implications and finally outline some new relations of the subject to algebraic
geometry.

\vskip 4.5mm


\end{abstract}

\vskip 12mm

\section{Introduction} \label{section 1}\setzero
\vskip-5mm \hspace{5mm }

Symplectic geometry exhibits a range of intersection phenomena
that cannot be predicted nor explained on the level of pure
topology or differential geometry. The main players in this game
are certain pairs of subspaces (e.g. Lagrangian submanifolds,
domains, or a mixture of both) whose mutual intersections cannot
be removed (or reduced) via the group of Hamiltonian or symplectic
diffeomorphisms.  The very first examples of such phenomena were
conjectures by Arnold in the 1960's, and eventually established
and further explored by Gromov, Floer and others starting from the
mid 1980s.

The first part of the paper will survey several intersection
phenomena and the mathematical tools leading to their discovery.
We shall not attempt to present the most general results and since
the literature is vast the exposition will be far from complete.
Rather we shall concentrate on various intersection phenomena
trying to understand their nature and whether there is any
relations between them.

The second part is dedicated to ``non-intersections'', namely to
situations where the principles of symplectic intersections break
down.  In the case of Lagrangian submanifolds this absence of
intersections is reflected in the vanishing of a symplectic
invariant called Floer homology. This vanishing when interpreted
algebraically leads to restrictions on the topology of Lagrangian
submanifolds. As a byproduct we shall explain how these
restrictions can be used to study some problems in algebraic
geometry concerning hyperplane sections and degenerations.

\section{Various intersection phenomena}\label{section 2}\setzero
\vskip-5mm \hspace{5mm }

In this section we shall make a brief tour through the zoo of
symplectic intersections, encountering three different species.

Before we start let us recall two important notions from
symplectic geometry. Let $(M, \omega)$ be a symplectic manifold. A
submanifold $L \subset M$ is called {\em Lagrangian} if $\dim L =
\frac{1}{2} \dim M$ and $\omega$ vanishes on $T(L)$. From now on
we assume all Lagrangian submanifolds to be closed. The second
notion is of {\em Hamiltonian
  isotopies}.  An isotopy of diffeomorphisms $\{h_t: M \to M\}_{0\leq
  t\leq 1}$, starting with $h_0 = \textnormal{id}$ is called
Hamiltonian if the (time-dependent) vector field $\xi_t$
generating it satisfies that the $1$-forms $i_{\xi_t} \omega$ are
{\em exact} for all $0\leq t \leq 1$.  Note that Hamiltonian
isotopies preserve the symplectic structure: $h_t^*\omega =
\omega$ for all $t$. Finally, two subsets $A, B \subset M$ are
said to be Hamiltonianly isotopic if there exists a Hamiltonian
isotopy $h_t$ such that $h_1(A)=B$. We refer the reader
to~\cite{M-S} for the foundations of symplectic geometry.

\subsection{Lagrangians intersect Lagrangians} \label{subsection 2.1}
\vskip-5mm \hspace{5mm }

 The most fundamental {\it Lagrangian
intersection} phenomenon occurs in cotangent bundles. Let $X$ be a
closed manifold and $T^*(X)$ be its cotangent bundle endowed with
the canonical symplectic structure $\omega_{\textnormal{can}}=\sum
dp_i \wedge dq_i$. Denote by $\lambda_{\textnormal{can}}=\sum p_i
dq_i$ the Liouville form (so that $\omega_{\textnormal{can}} =
d\lambda_{\textnormal{can}}$). Recall that a Lagrangian
submanifold $L \subset T^*(X)$ is called exact if the restriction
$\lambda_{\textnormal{can}}|_{T(L)}$ is exact. Note that the
property of exactness is preserved by Hamiltonian isotopies.
Denote by $O_X \subset T^*(X)$ the zero-section.  The following
theorem was proved by Gromov in~\cite{Gr}:

\noindent {\bf Theorem A.} \it Let $L \subset T^*(X)$ be an exact
Lagrangian submanifold. Then: \\
1) For every Lagrangian $L'$ which is Hamiltonianly isotopic to
$L$ we
have $L \cap L' \neq \emptyset$. \ \\
2) $L \cap O_X \neq \emptyset$. In particular, $L$ cannot be
separated from the zero-section by any Hamiltonian isotopy.  \rm

If one assumes $L$ to be a Hamiltonian image of the zero-section a
more quantitative version of Theorem~A holds:

\noindent {\bf Theorem B.} \it Let $L \subset T^*(X)$ be a
Lagrangian submanifold which is Hamiltonianly isotopic to the
zero-section and intersects it transversely. Then $$\# \, L \cap
O_X \geq \sum_{j=0}^{\dim X} b_j(X),$$ where $b_j(X)$ are the
Betti numbers of $X$. \rm

Chronologically Theorem~B preceded Theorem~A. It was conjectured
by Arnold (see~\cite{Ar-1} for the history), first proved for
$X=\mathbb{T}^n$ by Chaperon~\cite{Ch} and generalized to all
cotangent bundles by Hofer~\cite{Ho} and by Laudenbach and
Sikorav~\cite{Lau-Si}. Now a days it can be viewed as a special
case of Floer theory (see Section~2.4 below).

Note that the intersections described by both theorems above
cannot in general be understood on a purely topological level.
Indeed, in general topology predicts less than $\sum b_j(X)$
intersection points, and sometimes even none. Finally, note that
in general the statement of Theorem~B fails if one assumes $L$ to
be only symplectically isotopic to $O_X$, as the example
$X=\mathbb{T}^n$ shows.

\subsection{Balls intersect balls} \label{subsection 2.2}\vskip-5mm \hspace{5mm }

Denote by $B^{2n}(R)$ the closed Euclidean ball of radius $R$,
endowed with the standard symplectic structure induced from
$\mathbb{R}^{2n}$. Denote by ${\mathbb{C}}P^n$ the complex
projective space, endowed with its standard K\"{a}hler form
$\sigma$, normalized so that $\int_{\mathbb{C}P^1} \sigma = \pi$.
The following obstruction for {\it symplectic packing} was
discovered by Gromov~\cite{Gr}:

\noindent {\bf Theorem C.} \it Let $M$ be either $B^{2n}(1)$ or
${\mathbb{C}}P^n$. Let $B_{\varphi_1}, B_{\varphi_2} \subset M$ be
the images of two symplectic embeddings $\varphi_1:B^{2n}(R_1) \to
M$, $\varphi_2:B^{2n}(R_2) \to M$. If $R_1^2 + R_2^2 \geq 1$ then
$B_{\varphi_1} \cap B_{\varphi_2} \neq \emptyset$. \rm

Since symplectic embeddings are also volume preserving there is an
obvious volume obstruction for having $B_{\varphi_1} \cap
B_{\varphi_2} = \emptyset$. However, volume considerations predict
an intersection only if $R_1^{2n} + R_2^{2n} \geq 1$ (moreover for
volume preserving embeddings the latter inequality is sharp).

When one considers embeddings of several balls things become more
complicated and interesting. Here results are currently available
only in dimension $4$.

\noindent {\bf Theorem D.} \it Let $M$ be either $B^4(1)$ or
${\mathbb{C}}P^2$, and let $B_{\varphi_1}, \ldots, B_{\varphi_N}
\subset M$ be the images of symplectic embeddings
$\varphi_k:B^4(R) \to M$, $k=1,\ldots,N$, of $N$ balls of the same
radius $R$. Then there exist $i \neq j$ such that $B_{\varphi_i}
\cap B_{\varphi_j} \neq \emptyset$ in each of the following cases:
\begin{enumerate}
  \item $N=2$ or $3$ and $R^2 \geq 1/2$.
  \item $N=5$ or $6$ and $R^2 \geq 2/5$.
  \item $N=7$ and $R^2 \geq 3/8$.
  \item $N=8$ and $R^2 \geq 6/17$.
\end{enumerate}
Moreover all the above inequalities are sharp in the sense that in
each case if the inequality on $R$ is not satisfied then there
exist symplectic embeddings $\varphi_1, \ldots, \varphi_N$ as
above with disjoint images $B_{\varphi_1}, \ldots, B_{\varphi_N}
\subset M$. \rm

Statement 2 for $N=5$ was proved by Gromov~\cite{Gr}. The rest was
established by McDuff and Polterovich~\cite{M-P}. Let us mention
that for $N=4$ and any $N\geq 9$ this intersection phenomenon
completely disappears in the sense that an arbitrarily large
portion of the volume of $M$ can be filled by a disjoint union of
$N$ equal balls (see~\cite{M-P} for $N=4$ and $N=k^2$,
and~\cite{Bi-Pack, Bi-Stab} for the remaining cases).

\subsection{Balls intersect Lagrangians} \label{subsection 2.3}\vskip-5mm \hspace{5mm }

It turns out that there exist (symplectically) irremovable
intersections also between contractible domains (e.g. balls) and
Lagrangian submanifolds.

Denote by $\mathbb{R}P^n \subset {\mathbb{C}}P^n$ the Lagrangian
$n$-dimensional real projective space (embedded as the fixed point
set of the standard conjugation of ${\mathbb{C}}P^n$). The
following was proved in~\cite{Bi-Bar}:

\noindent {\bf Theorem E.} \it Let $B_{\varphi} \subset
{\mathbb{C}}P^n$ be the image of a symplectic embedding
$\varphi:B^{2n}(R) \to {\mathbb{C}}P^n$. If $R^2 \geq 1/2$ then
$B_{\varphi} \cap \mathbb{R}P^n \neq \emptyset$. Moreover the
inequality is sharp, namely for every $R^2 < 1/2$ there exists a
symplectic embedding $\varphi:B^{2n}(R) \to {\mathbb{C}}P^n$ whose
image avoids $\mathbb{R}P^n$.  \rm

In fact this pattern of intersections occurs in a wide class of
examples (see~\cite{Bi-Bar}):

\noindent {\bf Theorem E'.} \it Let $(M,\omega)$ be a closed
K\"{a}hler manifold with $[\omega] \in H^2(M;\mathbb{Q})$ and
$\pi_2(M)=0$. Then for every $\epsilon > 0$ there exists a
Lagrangian CW-complex $\Delta_{\epsilon} \subset (M, \omega)$ with
the following property: every symplectic embedding
$\varphi:B^{2n}(\epsilon) \to (M,\omega)$ must satisfy
$\textnormal{Image\,}(\varphi) \cap \Delta_{\epsilon} \neq
\emptyset$. \rm

By a Lagrangian CW-complex we mean a subspace $\Delta_{\epsilon}
\subset M$ which topologically is a CW-complex and the interior of
each of its cells is a smoothly embedded disc of $M$ on which
$\omega$ vanishes.

\subsection{Methods for studying intersections} 
{\bf Lagrangian intersections.} The first systematic study of
Lagrangian intersections was based on the theory of generating
function~\cite{Ch, Lau-Si} (an equivalent theory was independently
developed in contact geometry~\cite{Che}).  Gromov's theory of
pseudo-holomorphic curves~\cite{Gr} gave rise to an alternative
approach which culminated in what is now called Floer theory. Each
of these theories has its own advantage. Floer theory works in
larger generality and seems to have a richer algebraic structure,
on the other hand the theory of generating functions leads in some
cases to sharper results (see~\cite{E-G-Lag}).

Since Floer theory will appear in the sequel, let us outline a few facts about it (the reader is referred to the
works of Floer~\cite{F} and of Oh~\cite{Oh-Spectral, Oh-Monotone} for details). Let $(M, \omega)$ be a symplectic
manifold and $L_0, L_1 \subset (M, \omega)$ two Lagrangian submanifolds. In ``ideal'' situations Floer theory
assigns to this data an invariant $HF(L_0, L_1)$. This is a $\mathbb{Z}_2$-vector space obtained through an
infinite dimensional version of Morse-Novikov homology performed on the space of paths connecting $L_0$ to $L_1$.
The result of this theory is a chain complex $CF(L_0,L_1)$ whose underlying vector space is generated by the
intersection points $L_0 \cap L_1$ (one perturbs $L_0, L_1$ so their intersection becomes transverse). The
homology of this complex $HF(L_0,L_1)$ is called the Floer homology of the pair $(L_0, L_1)$. The most important
feature of $HF(L_0, L_1)$ is its invariance under Hamiltonian isotopies: if $L_0', L_1'$ are Hamiltonianly
isotopic to $L_0, L_1$ respectively, then $HF(L_0', L_1') \cong HF(L_0, L_1)$. From this point of view
$HF(L_0,L_1)$ can be regarded as a quantitative obstruction for Hamiltonianly separating $L_0$ from $L_1$. Indeed,
the rank of $HF(L_0, L_1)$ is a lower bound on the number of intersection points of any pair of transversally
intersecting Lagrangians $L_0',L_1'$ in the Hamiltonian deformation classes of $L_0, L_1$ respectively.

Let us explain the ``ideal situations'' in which Floer homology is
defined. First of all there are restrictions on $M$: due to
analytic difficulties manifolds are required to be either closed
or to have symplectically convex ends (e.g. $\mathbb{C}^n$,
cotangent bundles or any Stein manifold).  More serious
restrictions are posed on the Lagrangians. For simplicity we
describe them only for the case when $L_1$ is Hamiltonianly
isotopic to $L_0$. From now on we shall write $L=L_0$ and
$L'=L_1$. In Floer's original setting~\cite{F} the theory was
defined under the assumption that the homomorphism
$A_{\omega}:\pi_2(M,L) \to \mathbb{R}$, defined by $D \mapsto
\int_D \omega$, vanishes. The reason for this comes from the
construction of the differential of the Floer complex: the main
obstruction for defining a meaningful differential turns out to be
existence of holomorphic discs with boundary on $L$ or $L'$. These
discs appear as a source of non-compactness of the space of
solutions of the PDEs involved in the construction. Since
holomorphic discs must have positive symplectic area the
assumption $A_{\omega}=0$ rules out their existence. Under this
assumption Floer defined $HF(L,L')$ and proved its invariance
under Hamiltonian isotopies. Moreover he showed that $HF(L,L)$ is
isomorphic to the singular cohomology $H^*(L;\mathbb{Z}_2)$ of
$L$.  This together with the invariance give:

\noindent {\bf Theorem F.} \it Let $(M,\omega)$ be a symplectic
manifold, either compact or with symplectically convex ends. Let
$L \subset (M, \omega)$ be a Lagrangian submanifold with
$A_{\omega}=0$. Then for every Lagrangian $L'$ which is
Hamiltonianly isotopic to $L$ and intersects $L$ transversally we
have: $\# \, L \cap L' \geq \textnormal{rank} HF(L,L') =
\textnormal{rank} H^*(L; \mathbb{Z}_2)$.  In particular $L$ cannot
be separated from itself by a Hamiltonian isotopy. \rm

Floer theory was extended by Oh~\cite{Oh-Monotone} to cases when
$A_{\omega} \neq 0$. There are two assumptions needed for this
extension to work: the Maslov homomorphism $\mu:\pi_2(M,L) \to
\mathbb{Z}$ should be positively proportional to $A_{\omega}$
(such Lagrangians are called monotone). The second assumption is
that the positive generator $N_L$ of the subgroup
$\textnormal{Image\,}\mu \subset \mathbb{Z}$ is at least $2$. In
this setting Oh defined $HF(L,L')$ and proved its invariance under
Hamiltonian isotopies. It is however no longer true in general
that $HF(L,L)$ is isomorphic to $H^*(L;\mathbb{Z}_2)$. Still, Oh
proved~\cite{Oh-Spectral} that $HF(L,L)$ is related to
$H^*(L;\mathbb{Z}_2)$ through a spectral sequence. Recently the
theory was considerably generalized by Fukaya, Oh, Ohta and
Ono~\cite{FO3}.

\ \\{\bf Intersections of balls.}  Theorems~C and~D were obtained
using Gromov's theory of pseudo-holomorphic curves. The hard-core
of the proofs consists of existence of pseudo-holomorphic curves
of specified degrees that pass through a prescribed number of
points in the manifold (see~\cite{Gr,M-P}) for the details). From
a more modern perspective it can be viewed as an early application
of Gromov-Witten invariants.

Finally, Theorems~E and~E' are proved by a decomposition technique
introduced in~\cite{Bi-Bar} which enables to decompose symplectic
manifolds as a disjoint union of a symplectic disc bundle and a
Lagrangian $CW$-complex. A variation on the proof of Gromov's
non-squeezing theorem~\cite{Gr} gives an upper bound on the radius
of a symplectic ball that can be squeezed inside that disc bundle.
Hence, a larger ball must always intersect this $CW$-complex.  For
$M={\mathbb{C}}P^n$, the corresponding $CW$-complex turns out to
be a smooth copy of $\mathbb{R}P^n$. See~\cite{Bi-Bar} for the
details.

\section{Some questions and speculations}\label{section 3}\setzero
\vskip-5mm \hspace{5mm }

\noindent {\bf Cotangent bundles.} The following questions show
that even in the case of cotangent bundles the most fundamental
invariants are not completely understood.
\begin{enumerate}
  \item Let $L \subset T^*(X)$ be an exact Lagrangian (not necessarily
   Hamiltonianly isotopic to $O_X$). By Theorem~A, $L \cap O_X \neq
   \emptyset$. Is it true that $HF(L, O_X) \neq 0$~?
  \item Let $L_0, L_1 \subset T^*(X)$ be two exact Lagrangians (again,
   not necessarily Hamiltonianly isotopic neither to $O_X$ nor to each
   other). Is it true that $L_0 \cap L_1 \neq \emptyset$~? Is it true
   that $HF(L_0, L_1) \neq 0$~?
\end{enumerate}
These questions are of a theoretical importance, since the zero
section and its Hamiltonian images are the only known examples of
exact Lagrangians in $T^*(X)$.

\ \\ \noindent {\bf Symplectic packing.}  Lack of tools (or new
ideas) prevent us from understanding symplectic packings in
dimension higher than $4$. The only packing obstructions known in
these dimensions are described in Theorem~C. Note that
${\mathbb{C}}P^n$ admits full packing by $N=k^n$ equal
balls~\cite{M-P}, but it is unclear what happens for other values
of $N$. In view of this and Theorem~C, the first unknown case (for
$n\geq 3$) is of $N=2^n+1$ equal balls.

The situation in dimension $4$ is only slightly better. Except of
${\mathbb{C}}P^2$ and a few other rational surfaces no packing
obstructions are known. It is known that for every symplectic
$4$-manifold $(M, \omega)$ with $[\omega] \in H^2(M; \mathbb{Q})$
packing obstruction (for equal balls) disappear once the number of
balls is large enough (see~\cite{Bi-Stab}), but nothing is known
when the number of balls is small. In fact even the case of one
ball is poorly understood (namely, what is the maximal radius of a
ball that can be symplectically embedded in $M$). The reason here
is that the methods yielding packing obstructions strongly rely on
the geometry of algebraic and pseudo-holomorphic curves in the
manifold. The problem is that most symplectic manifolds have very
few (or none at all) $J$-holomorphic curves for a generic choice
of the almost complex structure. Thus, even in dimension $4$ it is
unknown whether or not packing obstructions is a phenomenon
particular to a sporadic class of manifolds such as
${\mathbb{C}}P^2$.

\ \\
\noindent {\bf Is everything Lagrangian?}  Weinstein's famous
saying could be relevant for the intersection described in
Theorems~C,D and~E. In other words, it could be that these
intersections are in fact Lagrangian intersections under disguise.
To be more concrete, let $\frac{1}{2} < R^2 < \frac{1}{2^{1/n}}$
and consider a Lagrangian $L_R$ lying on the boundary $\partial
B^{2n}(R)$. Is it possible to Hamiltonianly separate $L_R$ from
itself inside $B^{2n}(1)$ ?

If we can find a Lagrangian $L_R$ for which the answer is negative
then this would strongly indicate that Theorem~C is in fact a
Lagrangian intersections result. Namely it would imply Theorem~C
for $R_1 = R_2$ under the additional assumption that $\varphi_1,
\varphi_2$ are symplectically isotopic. A good candidate for $L_R$
seems to be the split torus $\partial B^2(\sqrt{R/n}) \times
\cdots \times \partial B^2(\sqrt{R/n}) \subset \partial
B^{2n}(R)$, but one could try other Lagrangians as well.

Attempts to approach this question with traditional Floer homology
fail. The reason is that Floer homology is blind to sizes: both to
the ``size'' of the Lagrangian $L_R$ as well as to the ``size'' of
the domain in which we work $B^{2n}(1)$. Indeed it is easy to see
that $HF(L_R, L_R)$ whether computed inside $B^{2n}(1)$ or in
$\mathbb{R}^{2n}$ is the same, hence vanishes. The meaning of
``sizes'' can be made precise: the size of $L_R$ is encoded in its
Liouville class, and the size of $B^{2n}(1)$ could be encoded here
by the action spectrum of its boundary.

It would be interesting to try a mixture of symplectic field
theory~\cite{SFT} with Floer homology. This would require a
sophisticated counting of holomorphic discs with $k$ punctures
(for all $k\geq 0$), where the boundary of the discs goe to $L_R$
and the punctures to periodic orbits on $\partial B^{2n}(1)$.

It is interesting to note that when the radii of the balls are not
equal things become more complicated. Indeed suppose that $R_1^2 +
R_2^2 > 1$ and consider two Lagrangian submanifolds $L_{R_1},
L_{R_2}$ lying on the boundaries of the balls $B_{\varphi_1},
B_{\varphi_2}$. Then clearly $L_{R_1}$ and $L_{R_2}$ can be
disjoint even though the balls $B_{\varphi_1}, B_{\varphi_2}$ do
intersect (e.g. two concentric balls $B_{\varphi_1} \subset
B_{\varphi_2}$, where $R_1 < R_2$). It would be interesting to see
to which extent this mutual position can be detected on the level
of the Lagrangians $L_{R_1}$ and $L_{R_2}$ alone. Or, in more
pictorial (but less mathematical) terms, do the Lagrangians
$L_{R_1}$ and $L_{R_2}$ know that they lie one ``inside'' the
other ?

Returning to the case of equal balls, if the above plan is
feasible, it would be interesting to try similar approaches for
more than two balls as described in Theorem~D. A similar approach
could be tried in the situation of Theorem~E. Here one could
expect an irremovable intersection between a Lagrangian
submanifold $L_R \subset \partial B_{\varphi}$ and
$\mathbb{R}P^n$.

\ \\ \noindent {\bf Quantitative intersections.} In contrast to
the quantitative version of Lagrangian intersections given by
Theorems~B and~F, Theorems~C--E provide only existence of
intersections. Is it possible to measure the size of these
intersections ?

More concretely, consider two balls $B_{\varphi_1}, B_{\varphi_2}
\subset B^{2n}(1)$ with $R_1^2 + R_2^2 > 1$ but with $R_1^{2n} +
R_2^{2n} < 1$ (so that $\textnormal{Vol}(B_{\varphi_1}) +
\textnormal{Vol} (B_{\varphi_2}) < 1$). Is it possible to bound
from below the size of $B_{\varphi_1} \cap B_{\varphi_2}$ ?

It is not hard to see that volume is a wrong candidate for the
size since for every $\epsilon > 0$ there exist two such balls
with $\textnormal{Vol}(B_{\varphi_1} \cap B_{\varphi_2}) <
\epsilon$. Symplectic capacities seem also to be inappropriate for
this task. It could be that ``size'' should be replaced here by a
kind of ``complexity'' or a trade-off between capacity and
complexity: namely if the intersection has large capacity (e.g.
when $B_{\varphi_1} \subset B_{\varphi_2}$) the complexity is low,
and vice-versa. Note that in dimension $2$ a possible notion of
complexity of a set is the number of connected components of its
interior.

A related problem is the following. Consider two symplectic balls
$B_{\varphi_1}, B_{\varphi_2} \subset {\mathbb{C}}P^n$ of radii
$R_1, R_2$, where $R_1^2 + R_2^2 = 1$. Assume further that
$\textnormal{Int\,}(B_{\varphi_1}) \cap
\textnormal{Int\,}(B_{\varphi_2}) = \emptyset$. Theorem~C implies
that the balls must intersect hence the intersection occurs on the
boundaries: $\partial B_{\varphi_1} \cap \partial B_{\varphi_2}
\neq \emptyset$. What can be said about the intersection $\partial
B_{\varphi_1} \cap \partial B_{\varphi_2} \neq \emptyset$, in
terms of size, dynamical properties etc. ?

It is easy to see that this intersection cannot be discrete.
Moreover, an argument based on the work of Sullivan~\cite{Su}
shows that the intersection must contain at least one entire
(closed) orbit of the characteristic foliation of the boundaries
of the balls (see~\cite{Po-Si} for a discussion on this point).
Looking at examples however suggests that the number of orbits in
the intersection should be much larger.

The same problem can be considered also for (some of) the extremal
cases described in Theorem~D. Similarly one can study the
intersection $\partial B_{\varphi} \cap \mathbb{R}P^n$ where
$B_{\varphi} \subset {\mathbb{C}}P^n$ is a symplectic ball of
radius $R^2=1/2$ whose interior is disjoint from $\mathbb{R}P^n$.
It is likely that methods of symplectic field theory~\cite{SFT}
could shed some light on this circle of problems.

\ \\ \noindent {\bf Stable intersections.}  The problems described
here come from Polterovich~\cite{Po-Diam}. Let $(M,\omega)$ be a
symplectic manifold and $A \subset M$ a subset. We say that $A$
has the {\it Hamiltonian intersection property} if for every
Hamiltonian diffeomorphism $f$ we have $f(A) \cap A \neq
\emptyset$. We say that $A$ has the {\it stable Hamiltonian
intersection property} if $O_{S^1} \times A \subset T^*(S^1)
\times M$ has the Hamiltonian intersection property. Polterovich
discovered in~\cite{Po-Diam} that if there exists a subset $A
\subset M$ with open non-empty complement and with the stable
Hamiltonian property then the universal cover
$\widetilde{\textnormal{Ham}}(M,\omega)$ of the group of
Hamiltonian diffeomorphisms has infinite diameter with respect to
Hofer's metric. Note that when $\pi_1(\textnormal{Ham}(M,\omega))$
is finite the same holds also for $\textnormal{Ham}(M,\omega)$
itself. (See~\cite{Po-Diam} for the details and references for
other results on the diameter of $\textnormal{Ham}$). This is
applicable when $(M,\omega)$ contains a Lagrangian submanifold $A$
with $HF(A,A) \neq 0$, since then $HF(O_{S^1}\times A,
O_{S^1}\times A) = (\mathbb{Z}_2 \oplus \mathbb{Z}_2) \otimes
HF(A,A) \neq 0$. For example, taking $A=\mathbb{R}P^n \subset
{\mathbb{C}}P^n$ Polterovich proved that $\textnormal{diam} \,
\widetilde{\textnormal{Ham}}({\mathbb{C}}P^n)=\infty$ (for $n=1,2$
the same holds for $diam \, {\textnormal{Ham}({\mathbb{C}}P^n)}$).

In view of the above the following question seems natural: does
every closed symplectic manifold contain a subset $A$ with open
non-empty complement and with the stable Hamiltonian intersection
property ? Note that besides Lagrangian submanifolds (with $HF
\neq 0$) no other stable Hamiltonian intersection phenomena are
known. It would also be interesting to find out whether the
intersections described in Theorems~C,D,E and especially $E'$
continue to hold after stabilization.

\section{Intersections versus non-intersections}\label{section 4}\setzero
\vskip-5mm \hspace{5mm }

In contrast to cotangent bundles there are manifolds in which {\it
  every compact subset} can be separated from itself by a Hamiltonian
isotopy. The simplest example is $\mathbb{C}^n$: indeed linear
translations are Hamiltonian, and any compact subset can be
translated away from itself. Clearly the same also holds for every
symplectic manifold of the type $M \times \mathbb{C}$ by applying
translations on the $\mathbb{C}$ factor. Note that manifolds of
the type $M \times \mathbb{C}$ sometime appear in ``disguised''
forms (e.g. as subcritical Stein manifolds, see
Cieliebak~\cite{Ci}).

The ``non-intersections'' property has quite strong consequences
on the topology of Lagrangian submanifolds already in
$\mathbb{C}^n$. Denote by $\omega_{\textnormal{std}}$ the standard
symplectic structure of $\mathbb{C}^n$ and let $\lambda$ be any
primitive of $\omega_{\textnormal{std}}$. Note that the
restriction $\lambda|_{T(L)}$ of $\lambda$ to any Lagrangian
submanifold $L \subset \mathbb{C}^n$ is closed. The following was
proved by Gromov in~\cite{Gr}:

\noindent {\bf Theorem G.} \it Let $L \subset \mathbb{C}^n$ be a
Lagrangian submanifold. Then the restriction of $\lambda$ to $L$
is not exact. In particular $H^1(L; \mathbb{R}) \neq 0$. \rm

Indeed if $\lambda$ were exact on $L$ then $A_{\omega}:
\pi_2(\mathbb{C}^n,L) \to \mathbb{R}$ must vanish, hence by
Theorem~F it is impossible to separate $L$ from itself by a
Hamiltonian isotopy. On the other hand, as discussed above, in
$\mathbb{C}^n$ this is always possible. We thus get a
contradiction.  (Gromov's original proof is somewhat different,
however a careful inspection shows it uses the failure of
Lagrangian intersections in an indirect way). Arguments exploiting
non-intersections were further used in clever ways by Lalonde and
Sikorav~\cite{La-Si} to obtain information on the topology of
exact Lagrangians in cotangent bundles (see also
Viterbo~\cite{Vi-Exact} for further results).

An important property of symplectic manifolds $W$ having the
``non-intersections'' property is the following vanishing
principle: {\it for every Lagrangian submanifold $L \subset W$
with well defined
  Floer homology we have $HF(L,L)=0$.} Applying this vanishing to
$\mathbb{C}^n$ yields restrictions on the possible Maslov class of
Lagrangian submanifolds of $\mathbb{C}^n$.  (Conjectures about the
Maslov class due to Audin appear already in~\cite{Au}. First
results in this directions are due to Polterovich~\cite{Po-Maslov}
and to Viterbo~\cite{Vi-Maslov}. The interpretation in
Floer-homological terms is due to Oh~\cite{Oh-Spectral}.
Generalizations to other manifolds appear in~\cite{ALP}
and~\cite{Bi-Ci-Stein}. Finally, consult~\cite{FO3} for recent
results answering old questions on the Maslov class).

\subsection{Lagrangian embeddings in closed manifolds}
\label{subsection 4.1} \vskip-5mm \hspace{5mm }

The ideas described above can be applied to obtain information on
the topology of Lagrangian submanifolds of some closed manifolds.
Note that in comparison to closed manifolds the case of
$\mathbb{C}^n$ can be regarded as local (Darboux Theorem). Of
course, ``local'' should by no means be interpreted as easy. On
the contrary, characterization of manifolds that admit Lagrangian
embeddings into $\mathbb{C}^n$ is completely out of reach with the
currently available tools.

Below we shall deal with the ``global'' case, namely with
Lagrangians in closed manifolds. One (coarse) way to ``mod out''
local Lagrangians is to restrict to Lagrangians $L$ with
$H_1(L;\mathbb{Z})$ zero or torsion (so that by Theorem~G they
cannot lie in a Darboux chart). The pattern arising in the
theorems below is that under such assumptions in some closed
symplectic manifolds we have homological uniqueness of Lagrangian
submanifolds. Let us view some examples.

We start with ${\mathbb{C}}P^n$. It is known that a Lagrangian
submanifold $L \subset {\mathbb{C}}P^n$ cannot have
$H_1(L;\mathbb{Z})=0$ (see Seidel~\cite{Se-Graded}, see
also~\cite{Bi-Ci-Closed} for an alternative proof). However,
$L\subset {\mathbb{C}}P^n$ may have torsion $H_1(L;\mathbb{Z})$ as
the example $\mathbb{R}P^n \subset {\mathbb{C}}P^n$ shows.

\noindent {\bf Theorem H.} \it Let $L \subset {\mathbb{C}}P^n$ be
a Lagrangian submanifold with $H_1(L;\mathbb{Z})$ a $2$-torsion
group (namely, $2 H_1(L;\mathbb{Z}) = 0$). Then:
\begin{enumerate}
  \item $H^*(L;\mathbb{Z}_2) \cong H^*(\mathbb{R}P^n; \mathbb{Z}_2)$
   as graded vector spaces.
  \item Let $a \in H^2({\mathbb{C}}P^n; \mathbb{Z}_2)$ be the
   generator. Then $a|_L \in H^2(L; \mathbb{Z}_2)$ generates the
   subalgebra $H^{\textnormal{even}}(L;\mathbb{Z}_2)$. Moreover if $n$
   is even the isomorphism in 1 is of graded algebras.
\end{enumerate}
\rm

Statement 1 of the theorem was first proved by
Seidel~\cite{Se-Graded}. An alternative proof based on
``non-intersections'' can be found in~\cite{Bi-Nonint}. Let us
outline the main ideas from~\cite{Bi-Nonint}.  Consider
${\mathbb{C}}P^n$ as a hypersurface of ${\mathbb{C}}P^{n+1}$. Let
$U$ be a small tubular neighbourhood of ${\mathbb{C}}P^n$ inside
${\mathbb{C}}P^{n+1}$. The boundary $\partial U$ looks like a
circle bundle over ${\mathbb{C}}P^n$ (in this case it is just the
Hopf fibration). Denote by $\Gamma_L \to L$ the restriction of
this circle bundle to $L \subset {\mathbb{C}}P^n$. A local
computation shows that $U$ can be chosen so that $\Gamma_L \subset
{\mathbb{C}}P^{n+1} \setminus {\mathbb{C}}P^n$ becomes a
Lagrangian submanifold. (This procedure works whenever we have a
symplectic manifold $\Sigma$ embedded as a hyperplane section in
some other symplectic manifolds $M$). The next observation is that
$\Gamma_L \subset {\mathbb{C}}P^{n+1} \setminus {\mathbb{C}}P^n$
is monotone and moreover its minimal Maslov number $N_{\Gamma_L}$
is the same as the one of $L$. Due to our assumptions on $H_1(L;
\mathbb{Z})$ this number turns out to satisfy $N_{\Gamma_L}\geq
n+1$. The crucial point now is that $HF(\Gamma_L, \Gamma_L)=0$.
Indeed, the symplectic manifold ${\mathbb{C}}P^{n+1} \setminus
{\mathbb{C}}P^n$ can be completed to be $\mathbb{C}^{n+1}$ where
Floer homology vanishes.

Having this vanishing we turn to an alternative computation of
$HF(\Gamma_L, \Gamma_L)$. This computation is based on the theory
developed by Oh~\cite{Oh-Spectral} for monotone Lagrangian
submanifolds.  According to~\cite{Oh-Spectral} Floer homology can
be computed via a spectral sequence whose first stage is the
singular cohomology of the Lagrangian. The minimal Maslov number
has an influence both on the grading as well as on the number of
steps it takes the sequence to converge to $HF$. In our case we
have a spectral sequence starting with $H^*(\Gamma_L;
\mathbb{Z}_2)$ and converging to $HF(\Gamma_L, \Gamma_L)=0$. A
computation through this process together with the information
that $N_{\Gamma_L} \geq n+1$ makes it possible to completely
recover $H^*(\Gamma_L; \mathbb{Z}_2)$. It turns out that
$H^i(\Gamma_L; \mathbb{Z}_2)=\mathbb{Z}_2$ for $i=0,1,n$ and
$n+1$, while $H^i(\Gamma_L; \mathbb{Z}_2)=0$ for all $1< i < n$.
Going back from $H^*(\Gamma_L; \mathbb{Z}_2)$ to $H^*(L;
\mathbb{Z}_2)$ is now done by the Gysin exact sequence of the
circle bundle $\Gamma_L \to L$ and noting that the second
Stiefel-Whitney class of this bundle is nothing but the
restriction $a|_L$ of the generator $a \in H^2({\mathbb{C}}P^n;
\mathbb{Z}_2)$.

Summarizing the proof, there are three main ingredients:
\begin{enumerate}
  \item Transforming the Lagrangian $L$ into a related Lagrangian
   $\Gamma_L$ living in a different manifold such that $\Gamma_L$ can
   be Hamiltonianly separated from itself. Consequently we obtain
   $HF(\Gamma_L, \Gamma_L)=0$.
  \item Relating $HF(\Gamma_L, \Gamma_L)$ to $H^*(\Gamma_L)$ via the
   theory of Floer homology (e.g. a spectral sequence).
  \item Passing back from $H^*(\Gamma_L)$ to $H^*(L)$.
\end{enumerate}

Similar ideas work in various other cases (see~\cite{Bi-Nonint}).
For example, consider ${\mathbb{C}}P^n \times {\mathbb{C}}P^n$.
This manifold has Lagrangians with $H_1(L;\mathbb{Z})=0$, e.g.
${\mathbb{C}}P^n$ which can be embedded as the ``anti-diagonal''
$\{(z,w) \in {\mathbb{C}}P^n \times {\mathbb{C}}P^n |
w=\overline{z} \}$.

\noindent {\bf Theorem I.} \it Let $L \subset {\mathbb{C}}P^n
\times {\mathbb{C}}P^n$ be a Lagrangian with $H_1(L;
\mathbb{Z})=0$. Then $H^*(L; \mathbb{Z}_2) \cong
H^*({\mathbb{C}}P^n; \mathbb{Z}_2)$, the isomorphism being of
graded algebras. \rm

Another application of this circle of ideas is for Lagrangian
spheres. Recently Lagrangian spheres have attracted special
attention due to their relations to interesting symplectic
automorphisms~\cite{Se-Thesis, Se-Graded} and to symplectic
Lefschetz pencils~\cite{Do-Polynom}.

\noindent {\bf Theorem J.} \it 1) Let $M$ be a closed symplectic
manifold with $\pi_2(M)=0$, and denote by $m=\dim_{\mathbb{C}}M$
its complex dimension. If $M \times {\mathbb{C}}P^n$ (where
$m,n\geq 1$) has a
Lagrangian sphere then $m \equiv n+1 \pmod{2n+2}$. \ \\
\noindent 2) Let $M={\mathbb{C}}P^n \times {\mathbb{C}}P^m$,
$m+n\geq 3$, be endowed with the split symplectic form
$(n+1)\sigma \oplus (m+1)\sigma$. If $M$ has a Lagrangian sphere
then $\gcd(n+1, m+1)=1$. \rm

Let us remark that when $m=n+1$ any product of the form
${\mathbb{C}}P^n \times M$ (with $\dim_{\mathbb{C}}M=n+1$) indeed
has a Lagrangian sphere, after a possible rescaling of the
symplectic form on the $M$ factor
(see~\cite{ALP},~\cite{Bi-Ci-Closed}). We are not aware of any
other examples, namely when $m \equiv n+1 \pmod{2n+2}$ but $m\neq
n+1$).

\section{Relations to algebraic geometry}\label{section 5}\setzero
\vskip-5mm \hspace{5mm }

The purpose of this section is to show how ideas from Section~4
are related to algebraic geometry. We shall not present new
results here but rather try to outline a new direction in which
symplectic methods can be used in algebraic geometry.

\subsection{Hyperplane sections} \label{subsection 5.1}\vskip-5mm \hspace{5mm }

Let $\Sigma$ be a smooth projective variety. The classical
Lefschetz theorem provides restrictions on smooth varieties $X$
that may contain $\Sigma$ as their hyperplane section. It was
discovered by Sommese~\cite{So} that there exist projective
varieties $\Sigma$ that cannot be hyperplane sections (or even
ample divisors) in {\it any} smooth variety $X$. For example,
Sommese proved that Abelian varieties of (complex) dimension $\geq
2$ have this property (see~\cite{So} for more examples).

Let us outline an alternative approach to this problem using
symplectic geometry. Let $X \subset {\mathbb{C}}P^N$ be a smooth
variety. Denote by $X^{\vee} \subset ({\mathbb{C}}P^N)^*$ the dual
variety (namely, the variety of all hyperplanes $H \in
({\mathbb{C}}P^N)^*$ that are non-transverse to $X$).

\noindent {\bf Theorem K.} \it Suppose that $\Sigma = X \cap H_0
\subset X$ is a smooth hyperplane section of $X$ obtained from a
projective embedding $X \subset {\mathbb{C}}P^N$. Then either
$\Sigma$ has a Lagrangian sphere (for the symplectic structure
induced from ${\mathbb{C}}P^N$), or
$\textnormal{codim}_{\mathbb{C}}(X^{\vee}) > 1$.  \rm

Here is an outline of the proof. Suppose that
$\textnormal{codim}_{\mathbb{C}}(X^{\vee}) = 1$. Choose a generic
line $\ell \subset ({\mathbb{C}}P^N)^*$ intersecting $X^{\vee}$
transversely (and only at smooth points of $X^{\vee}$). Consider
the pencil $\{X \cap H\}_{H \in \ell}$ parametrized by $\ell$.
Passing to the blow-up $\widetilde{X}$ of $X$ along the base locus
of the pencil we obtain a holomorphic map $\pi:\widetilde{X} \to
\ell \approx \mathbb{C}P^1$. The critical values of $\pi$ are in
1-1 correspondence with the point of $\ell \cap X^{\vee}$.
Moreover, the fact that $\ell$ intersects $X^{\vee}$ transversely
implies that $\pi$ is a so called Lefschetz fibration, namely each
critical point of $\pi$ has non-degenerate (complex) Hessian (in
other words, locally $\pi$ looks like a holomorphic Morse
function). The condition
$\textnormal{codim}_{\mathbb{C}}(X^{\vee}) = 1$ ensures that $\ell
\cap X^{\vee} \neq \emptyset$ hence at least one of the fibres of
$\pi$ is singular. Let $X_0$ be such a fibre and $p \in X_0$ a
critical point of $\pi$. The important point now is that the
vanishing cycle (corresponding to $p$) that lies in the nearby
smooth fibre $X_{\epsilon}$ can be represented by a (smooth)
Lagrangian sphere. By Moser argument all the smooth divisors in
the linear system $\{ X \cap H \}_{H \in ({\mathbb{C}}P^N)^*}$ are
symplectomorphic. In particular $\Sigma$ has a Lagrangian sphere
too.

The existence of Lagrangian vanishing cycles was known
folklorically for long time. Its importance to symplectic geometry
was realized by Arnold~\cite{Ar-2}, Donaldson~\cite{Do-Polynom}
and by Seidel~\cite{Se-Thesis}.

Theorem~K can be applied as follows: given a smooth variety
$\Sigma$, use methods of symplectic geometry to prove that
$\Sigma$ contains no Lagrangian spheres, say for any symplectic
structure compatible with the complex structure of $\Sigma$. Then
by Theorem~K the only chance for $\Sigma$ to be a hyperplane
section is inside a variety $X$ with ``small dual'', namely
$\textnormal{codim}_{\mathbb{C}}(X^{\vee})>1$. Let us remark that
smooth varieties $X\subset {\mathbb{C}}P^N$ with
$\textnormal{codim}_{\mathbb{C}}(X^{\vee}) > 1$ are quite rare,
and have very restricted geometry (see e.g. Zak~\cite{Z} and
Ein~\cite{Ei-1,Ei-2}). Using the theory of ``small dual
varieties'' we can either rule out this case or get strong
restrictions on the pair $(X,\Sigma)$.

Let us illustrate this on the example mentioned at the beginning
of the section. Let $\Sigma$ be an Abelian variety of complex
dimension $n\geq 2$. Note that $\Sigma$ cannot have a Lagrangian
sphere for any K\"{a}hler form. Indeed, if $\Sigma$ had such a
sphere then the same would hold also for the universal cover of
$\Sigma$ which is symplectomorphic to $\mathbb{C}^n$. But this is
impossible in view of Theorem~G. Thus if $\Sigma$ is a hyperplane
section of $X \subset {\mathbb{C}}P^N$ then
$\textnormal{codim}_{\mathbb{C}}(X^{\vee}) > 1$. It is well
known~\cite{Kl} that in this case $X$ must have rational curves
(in fact lots of them). In particular $\pi_2(X) \neq 0$. By
Lefschetz's theorem we get $\pi_2(\Sigma) \neq 0$. But this is
impossible since $\Sigma$ is an Abelian variety. We therefore
conclude that $\Sigma$ cannot be a hyperplane section in {\it any}
smooth variety $X$.

An analogous (though symplectically more involved) argument should
apply also to any algebraic variety $\Sigma$ with $c_1=0$ and
$b_1(\Sigma) \neq 0$ (see~\cite{Bi-Algebraic}). An application of
more refined symplectic tools (e.g. methods described in
Section~4.1 above) can be used to obtain many more examples.

Here is another typical application: let $C$ be a projective curve
of genus $>0$.  It was observed by Silva~\cite{Si} that $C \times
{\mathbb{C}}P^n$ can be realized as a hyperplane section in
various smooth varieties. Note that by Theorem~J, $C \times
{\mathbb{C}}P^n$ cannot have any Lagrangian spheres. It
immediately follows that the only smooth varieties $X$ that
support $C \times {\mathbb{C}}P^n$ as their hyperplane section
must have small dual. For $n\leq 5$ results of
Ein~\cite{Ei-1,Ei-2} make it even possible to list all such $X$'s.

We conclude with a remark on the methods. The symplectic approach
outlined above gives coarser results. Indeed Sommese~\cite{So}
provides examples of varieties that cannot be ample divisors
whereas the methods above only rule out the possibility of being
very ample. On the other hand the symplectic approach has an
advantage in its robustness with respect to small deformations
(see~\cite{Bi-Algebraic}, c.f.~\cite{So-Singular}).

\subsection{Degenerations of algebraic varieties} \label{subsection 5.2}
\vskip-5mm \hspace{5mm }

 The methods of the previous section can
also be used to study degenerations of algebraic varieties. Let
$Y$ be a smooth projective variety. We say that $Y$ admits a
K\"{a}hler degeneration with isolated singularities if there
exists a K\"{a}hler manifold $X$ and a proper holomorphic map
$\pi:X \to D$ to the unit disc $D \subset \mathbb{C}$ with the
following properties:
\begin{enumerate}
  \item Every $0 \neq t \in D$ is a regular value of $\pi$ (hence, all
   the fibres $X_t=\pi^{-1}(t)$, $t \neq 0$, are smooth K\"{a}hler
   manifolds).
  \item $0$ is a critical value of $\pi$ and all the critical points
   of $\pi$ are isolated.
  \item $Y$ is isomorphic (as a complex manifold) to one of the smooth
   fibres of $\pi$, say $X_{t_0}$, $t_0 \neq 0$.
\end{enumerate}

As in the previous section this situation is related to symplectic
geometry through the Lagrangian vanishing cycle construction. As
pointed out by Seidel~\cite{Se-Graded} one can locally morsify
each of the critical points in $X_0=\pi^{-1}(0)$ and then by
applying Moser's argument obtain for each critical point of $\pi$
at least one Lagrangian sphere in the nearby fibre $X_{\epsilon}$.
Since all the smooth fibres are symplectomorphic we obtain
Lagrangian spheres also in $Y$.

Applying results from Section~4 to this situation we obtain
examples of projective varieties that do not admit {\em any}
degeneration with isolated singularities. For example, let $Y$ be
any of the following:
\begin{itemize}
  \item ${\mathbb{C}}P^n$, $n \geq 2$. Or more generally
   ${\mathbb{C}}P^n \times M$, where $M$ is a smooth variety with
   $\pi_2(M)=0$ and $\dim_{\mathbb{C}}M \not \equiv n+1 \pmod{n+1}$.
  \item Any variety whose universal cover is $\mathbb{C}^n$, ($n\geq
   2$), or a domain in $\mathbb{C}^n$.
\end{itemize}
Then by the results in Section~4, $Y$ has no Lagrangian spheres,
hence does not admit any degeneration as above. More examples can
be found in~\cite{Bi-Algebraic}.

This point of view seems non-trivial especially when
$H_n(Y;\mathbb{Z})=0$, where $n=\dim_{\mathbb{C}}Y$. In these
cases the vanishing cycles are zero in homology and it seems that
there are no obvious topological obstructions for degenerating $Y$
as above.
From the list above, the first non-trivial example should be
${\mathbb{C}}P^n$ with $n=\textnormal{odd}\geq 3$. It would be
interesting to figure out to which extent the above statement
could be proved within the tools of pure algebraic geometry. Note
that Lagrangian spheres are a non-algebraic object and it seems
that their existence/non-existence cannot be formalized in purely
algebro-geometric terms.

Another direction of applications should be to find an upper bound
on the number of singular points of an algebraic variety $X_0$
that can be obtained from a degeneration of $Y$.  Note that the
vanishing cycles of different singular points of $X_0$ are
disjoint. Thus the idea here is to obtain an upper bound on the
number of possible disjoint Lagrangian spheres that can be
embedded in $Y$. The simplest test case here should be the quadric
$Q=\{z_0^2 + \cdots + z_{n+1}^2 = 0\} \subset
{\mathbb{C}}P^{n+1}$, where $n\geq 2$. Clearly $Q$ can be
degenerated to a variety $X_0$ with isolated singularities (e.g.
to a cone over a smaller dimensional quadric). It seems reasonable
to expect that in every such degeneration the singular fibre $X_0$
will have only one singular point. Note that for
$n=\textnormal{even}$ this easily follows from topological reason
but it may not be so when $n=\textnormal{odd}\geq 3$ because
$H_n(Q; \mathbb{Z})=0$.  From a symplectic point of view the above
statement would follow if we could prove that every two Lagrangian
spheres in $Q$ must intersect.  This is currently still unknown
but there are evidences supporting this
conjecture~\cite{Bi-Nonint}. It is likely that a refinement of the
methods from~\cite{Se-Long} would be useful for this purpose.
More generally, one could try to bound the number of singular
fibres in a degeneration of other hypersurfaces $\Sigma \subset
{\mathbb{C}}P^{n+1}$ (in terms of $\deg(\Sigma)$ and $n$).
See~\cite{Bi-Nonint, Bi-Algebraic} for the conjectured bounds.

\ \\{\bf Acknowledgments.}  I am indebted to Leonid Polterovich
for sharing with me his insight into symplectic intersections in
general and especially regarding the problems presented in
Section~3. Numerous discussions with him have influenced my
conceptions on the subject. I wish to thank him also for valuable
comments on earlier drafts of the paper. I would like to thank
also Jonathan Wahl who told me about Sommese's work~\cite{So},
Olivier Debarre for the reference to the works of Zak~\cite{Z} and
Ein~\cite{Ei-1,Ei-2}, and Paul Seidel for useful discussions on
symplectic aspects of singularity theory and the material
presented in section~5.

\label{lastpage}

\end{document}